\documentclass{amsart}

\newcommand{\nt}{\noindent}

\newcommand{\bC}{{\mathbb C}}
\newcommand{\bP}{{\mathbb P}}
\newcommand{\bQ}{{\mathbb Q}}
\newcommand{\bR}{{\mathbb R}}
\newcommand{\bZ}{{\mathbb Z}}
\newcommand{\Ext}{\mbox{Ext}}
\newcommand{\cI}{{\mathcal I}}
\newcommand{\cO}{{\mathcal O}}
\newcommand{\IZ}{{\cI}_Z}
\newcommand{\OS}{{\mathcal O}_S}

\newcommand{\DX}{{\mathcal D}(X)}
\newcommand{\IW}{{\mathcal I}_W}
\newcommand{\IC}{{\mathcal I}_C}
\newcommand{\OC}{{\mathcal O}_C}
\newcommand{\IV}{{\mathcal I}_V}
\newcommand{\OX}{{\mathcal O}_X}
\newcommand{\rk}{{\mbox{rk}}}
\newcommand{\F}{{\mathcal F}}
\newcommand{\T}{{\mathcal T}}
\newcommand{\A}{{\mathcal A}}
\newcommand{\D}{{\mathcal D}}

\begin{document}

\title{Reider's Theorem and Thaddeus Pairs Revisited}

\author{Daniele Arcara}

\address{Department of Mathematics, Saint Vincent College,
300 Fraser Purchase Road, Latrobe, PA 15650-2690, USA}

\email{daniele.arcara@email.stvincent.edu}

\author{Aaron Bertram}

\address{Department of Mathematics, University of Utah,
155 S. 1400 E., Room 233, Salt Lake City, UT 84112-0090, USA}

\email{bertram@math.utah.edu}

\thanks{The second author was partially supported by NSF grant
DMS 0501000}

\maketitle


\nt {\bf 1. Introduction.} Let $X$ be a smooth projective variety over $\bC$ of dimension $n$ equipped
with an ample line bundle $L$ and a subscheme $Z \subset X$ of length $d$.
Serre duality provides a natural isomorphism of vector spaces (for each $i = 0,...,n$):
$$(*) \ \ \Ext^i(L\otimes \IZ, \OX)  \cong \mbox{H}^{n-i}(X,K_X\otimes L \otimes \IZ)^\vee$$

Thaddeus pairs and Reider's theorem concern the cases $i=1$  and $n= 1, 2$.
In these cases one associates a rank two coherent sheaf $E_\epsilon$
to each {\it extension class} $\epsilon \in \Ext^1(L\otimes \IZ, \OX)$ via the  short exact sequence:
$$(**) \ \ \epsilon: 0 \rightarrow \OX \rightarrow 
E_\epsilon \rightarrow L\otimes \IZ \rightarrow 0$$
and the Mumford stability (or instability) of $E_\epsilon$ allows one  
to distinguish among extension classes. The ultimate aim of this paper is to show how
a new notion of Bridgeland stability can similarly be used to distinguish among higher extension classes,
leading to a natural higher-dimensional 
generalization of Thaddeus pairs as well as the setup for a higher-dimensional Reider's theorem.

\medskip

Reider's theorem gives numerical conditions on an ample  line bundle $L$ on 
a surface $S$ that guarantee the vanishing of the vector spaces
$\mbox{H}^1(S,K_S \otimes L \otimes \IZ)$
which in turn implies the base-point-freeness (the $d=1$ case) and very ampleness 
(the $d=2$ case) 
of the adjoint line bundle $K_S \otimes L$.

\medskip

In the first part of this note we will revisit Reider's Theorem in the context of Bridgeland stability conditions. 
Reider's approach, following Mumford,
uses the Bogomolov inequality for Mumford-stable coherent sheaves on a surface to argue  
(under suitable numerical conditions on $L$) that
no exact sequence $(**)$ can produce a Mumford stable 
sheaf $E_\epsilon$, and then uses the Hodge Index Theorem to argue that the only 
exact sequences $(**)$ that produce non-stable sheaves must split.
Thus one concludes that $\mbox{Ext}^1(L\otimes \IZ,\OS) = 0$ and $\mbox{H}^1(S,K_S \otimes L \otimes \IZ) = 0$,
as desired.

\medskip

Here, we will regard an extension class in $(**)$  as a morphism to the {\it shift} of $\OS$:
$$\epsilon: L\otimes \IZ \rightarrow \OS[1]$$
in one of a family of {\it tilts} $\A_s \ (0 < s < 1)$ of the abelian category of coherent sheaves on $X$ within the 
bounded derived category $\DX$ of complexes of coherent sheaves on $X$. Reider's argument for a surface $S$ is essentially equivalent to 
ruling out non-trivial extensions by determining that:

\medskip

$\bullet$ $\epsilon$ is neither injective nor surjective and

\medskip

$\bullet$ if neither injective nor surjective, then $\epsilon = 0$ (using Hodge Index).

\medskip

This way of looking at Reider's argument allows for some minor improvements, but more importantly leads to the notion of  {\it Bridgeland stablility conditions}, which are stability 
conditions, not on coherent sheaves, but rather on objects of $\A_s$.

\medskip

In \cite{AB}, it was shown that the 
Bogomolov Inequality and Hodge Index Theorem imply the existence of such stability conditions on
arbitrary smooth projective surfaces $S$ (generalizing Bridgeland's stability conditions for $K$-trivial surfaces \cite{Bri08}). Using these
stability conditions, we investigate the stability of objects of the form
$L\otimes \IZ$ and $\OS[1]$
with a view toward reinterpreting the vanishing:
$$\mbox{Hom}(L\otimes\IZ , \OS[1]) = 0$$  
as a consequence of an inequality $\mu(L\otimes \IZ) > \mu(\OS[1])$ of 
Bridgeland slopes. Since this is evidently a stronger 
condition than just the vanishing of the Hom, it is unsurprising that it should require stronger
numerical conditions. This reasoning easily generalizes to the case where $\OS$ is replaced
by $\IW^\vee$, the derived dual of the ideal sheaf of a finite length subscheme $W \subset S$.

\medskip

The Bridgeland stability of the objects $L\otimes \IZ$ and $\IW^\vee[1]$ is central to 
a new generalization of Thaddeus pairs from curves to surfaces. 
A Thaddeus pair on a curve $C$ is an extension of the form:
$$\epsilon: 0 \rightarrow \IW^\vee \rightarrow E_\epsilon \rightarrow L\otimes \IZ \rightarrow 0$$
where $L$ is a line bundle and $Z,W \subset C$ are effective divisors. 
Normally we would write this:
$$\epsilon: 0 \rightarrow \OC(W) \rightarrow E_\epsilon \rightarrow L(-Z) \rightarrow 0$$
since finite length subschemes of a curve are effective Cartier divisors.
The {\it generic} such extension determines a Mumford-stable 
vector bundle $E_\epsilon$ on $C$ whenever:
$$\deg(L(-Z)) > \deg(\OC(W)) \ \ (\mbox{and} \ \ C \ne {\bP}^1)$$
or, equivalently, whenever the Mumford slope of $L(-Z)$ exceeds that of $\OC(W)$ 
(both line bundles are trivially Mumford-stable). Moreover, the Mumford-unstable 
vector bundles arising in this way are easily described
in terms of the secant varieties to the image of $C$ under the natural linear series map:
$$\phi: C \rightarrow {\bP}(\mbox{H}^0(C,K_C\otimes L(-Z-W))^\vee)
\cong {\bP}(\Ext^1(L(-Z),\OC(W))) $$
since an unstable vector bundle $E_\epsilon$  can only be destabilized 
by a sub-line bundle $L(-Z') \subset L(-Z)$ that lifts to a sub-bundle of $E_\epsilon$:
$$\begin{array}{ccccccccc}
&&&&&& L(-Z')\\
&&&&&\swarrow & \downarrow\\
(\dag) \ \ \  0 & \rightarrow & \OC(W) & \rightarrow & E_\epsilon & \rightarrow & L(-Z) & \rightarrow & 0
\end{array}$$

\medskip

In the second part of this paper, we note that Thaddeus pairs naturally generalize to surfaces as extensions of the form:
$$\epsilon: 0 \rightarrow \IW^\vee[1] \rightarrow E^\bullet _\epsilon
 \rightarrow L\otimes \IZ \rightarrow 0$$
in the categories $\A_s$ under appropriate Bridgeland stability conditions for which 
both $L\otimes \IZ$ and $\IW^\vee[1]$ are Bridgeland stable and their Bridgeland 
slopes satisfy
$\mu(L\otimes \IZ) > \mu(\IW^\vee[1]).$ Note that $E^\bullet_\epsilon$ 
is not ever a coherent sheaf. 

\medskip

This is a very satisfying generalization of Thaddeus pairs since:
$$\mbox{Ext}^1_{\A_s}(L\otimes \IZ, \IW^\vee[1]) \cong 
{\mbox H}^0(S,K_S \otimes L\otimes \IZ\otimes \IW)^\vee$$
by Serre duality. In this case, however, there are subobjects:
$$K \subset L\otimes\IZ$$
not of the form $L\otimes \cI_{Z'}$ that may destabilize $E^\bullet_\epsilon$, 
as in $(\dag)$. These subobjects are necessarily
coherent sheaves, but may be of higher rank than one, and therefore not subsheaves of
$L\otimes \IZ$ in the usual sense. This leads to a much richer geometry for the locus
of ``unstable'' extensions than in the curve case.

\medskip

We will finally discuss the moduli problem for families 
of Bridgeland stable objects with the particular invariants:
$$[E] = [L\otimes \IZ] + [\IW^\vee[1]] = [L\otimes \IZ] - [\IW^\vee]$$
in the Grothendieck group (or cohomology ring) of $S$, and finish by describing 
wall-crossing phenomena of (some of) these moduli spaces in the $K$-trivial case, following \cite{AB}.

\medskip

This line of reasoning suggests a natural question for {\it three-folds} $X$. Namely, might it 
be possible to prove a Reider theorem for $L$ and $Z \subset X$ by ruling out 
non-trivial extensions of 
the form:
$$\epsilon: 0 \rightarrow \OX[1] \rightarrow E^\bullet _\epsilon 
\rightarrow L\otimes \IZ \rightarrow 0$$
in some tilt $\A_s$ of the category of  coherent sheaves on $X$ via a version of the 
Bogomolov Inequality and Hodge Index Theorem for objects of $\A_s$ on threefolds? 

\medskip

We do not 
know versions of these results that would allow a direct application of 
Reider's method of proof,
but this seems a potentially fruitful direction for further research, and ought to be
related to the current active search for examples of  Bridgeland stability conditions on complex projective threefolds.

\bigskip

\nt {\bf 2. The Original Reider.}  Fix an ample divisor $H$ on a 
smooth projective variety $X$ over $\bC$ of dimension $n$. 
A torsion-free coherent sheaf $E$ on $X$ has Mumford slope:
$$\mu_H(E) = \frac{c_1(E)\cdot H^{n-1}}{\rk(E) H^n}$$
and $E$ is $H$-{\it Mumford-stable} if
$\mu_H(K) < \mu_H(E)$
for all subsheaves $K \subset E$ with the property that $Q = E/K$ is supported in codimension $\le 1$.

\medskip

\nt {\bf Bogomolov Inequality:} Suppose $E$ is $H$-Mumford-stable and $n \ge 2$. Then:
$$\mbox{ch}_2(E)\cdot H^{n-2}  \le \frac{c_1^2(E) \cdot H^{n-2}}{2\rk(E)}$$
(in  case $X = S$ is a surface, the conclusion is independent of the choice of $H$)

\medskip

\nt {\bf Application 2.1:} For $S,L,Z$ as above, suppose $\epsilon \in \mbox{Ext}^1(L\otimes \IZ,\OS)$ and:
$$\epsilon: 0 \rightarrow \OS \rightarrow E \rightarrow L\otimes \IZ \rightarrow 0$$
yields a Mumford-stable sheaf $E$. Then
$c_1^2(L) \le 4d$.

\medskip

{\bf Proof:} By the Bogomolov inequality:
$$\mbox{ch}_2(E) = \frac{c_1^2(L)}2 - d \le \frac{c_1^2(E)}{4} = \frac {c_1^2(L)}4$$

\nt {\bf Hodge Index Theorem:} Let $D$ be an arbitrary divisor on $X$. Then:
$$\left(D^2\cdot H^{n-2}\right)\left(H^n\right) \le (D\cdot H^{n-1})^2$$
and equality holds if and only if 
$D\cdot E\cdot H^{n-2}= k E\cdot H^{n-1}$ for some $k\in \bQ$ and all divisors $E$.

\medskip

\nt {\bf Application 2.2:} Suppose $E$ is not $c_1(L)$-Mumford-stable in Application 2.1. 
Then  either $\epsilon = 0$ or else there is an 
effective curve $C \subset S$ such that:
$$\mbox{(a)} \ \ \ C\cdot c_1(L) \le \frac 12 c_1^2(L) \ \ \ \ \mbox{and} \ \  \mbox{(b)}  \ \ \ 
C\cdot c_1(L) \le C^2 + d\ \ \ \ \ \ \ \ \ \ \ $$
and it follows that $-d < C^2 \le d$. Moreover, 
$$\mbox{(c)} \ \ c_1^2(L) > 4d \Rightarrow C^2 < d \ \ \mbox{and}\ \ 
\mbox{(d)} \ \ c_1^2(L) > (d+1)^2 \Rightarrow C^2 \le 0$$

\medskip

{\bf Proof:} By definition of (non)-stability, there is a 
rank-one subsheaf $K \subset E$ such that
$c_1(K)\cdot c_1(L) \ge \frac 12 c_1(E)\cdot c_1(L) = \frac 12 c_1^2(L)$. It follows that the induced map
$K \rightarrow L\otimes \IZ$ is non-zero and either $K$ splits the sequence, or else 
$K \subset L\otimes \IZ$ is a proper subsheaf. In the latter case, $K = L(-C)\otimes \IW$ for 
some effective curve $C$ and zero-dimensional $W \subset S$, and (a) now follows immediately.

\medskip

The inequality (b) is seen by computing the second Chern character of $E$ in two
different ways. We may assume
without loss of generality that the cokernel $Q= E/K$ is also torsion-free by replacing $Q$ with its 
torsion-free quotient $Q'$ and $K$ with the kernel of the induced map $E \rightarrow Q'$ if necessary
(this will only increase the value of $c_1(K)\cdot c_1(L)$). Then $Q$ has the form $\OS(C)\otimes \IV$ for 
some $V \subset S$, and in particular:
$$\mbox{ch}_2(E) = \frac {c_1^2(L)}2 - d = \frac{(c_1(L) - C)^2}2 - l(W) + \frac{C^2}2 - l(V) 
\le \frac{c_1^2(L)}2 + C^2 - C\cdot c_1(L)$$
which gives (b).

\medskip

Next, applying Hodge Index to (a) and (b) gives:
$$C^2c_1^2(L) \le (C\cdot c_1(L))^2 \le \frac 12 c_1^2(L)\left(C^2 + d\right)$$
from which we conclude that $C^2 \le d$. That $C^2 > -d$ follows immediately from (b) and 
the fact that $L$ is ample.
Finally, let $C^2 = d -k $ for $0 \le k < d$ and apply Hodge Index to (b) to conclude that:
$$(d-k)c_1^2(L) \le (C\cdot c_1(L))^2 \le (2d - k)^2 \Rightarrow c_1^2(L) \le 4d + \frac{k^2}{d-k}$$
and then (c) and (d) follow from the cases $k=0$ and $k \le d-1$, respectively.

\medskip

All of this gives as an immediate corollary a basic version of:

\medskip

\nt {\bf Reider's Theorem:} If $L$ is an ample line bundle
on a smooth projective surface $S$ such that $c_1^2(L) > (d+1)^2$ and $C\cdot c_1(L) > C^2 + d$ 
for all effective divisors $C$ on $S$ satisfying $C^2 \le 0$, then 
``$K_S + L$ separates length $d$ subschemes of $S$,'' i.e.
$$\mbox{H}^1(S, K_S\otimes L \otimes \IZ) = 0$$
for all subschemes $Z \subset S$ of length $d$ (or less). 

\medskip

\nt {\bf Corollary (Fujita's Conjecture for Surfaces):} If $L$ is an ample line bundle on a smooth projective surface $S$, then $K_S + (d+2)L$ separates length $d$ subschemes.

\medskip

\nt {\bf Note:} For other versions of Reider's theorem, see e.g.\ \cite{Laz97}. 

\bigskip

\nt {\bf 3. Reider Revisited.} A torsion-free coherent sheaf $E$ is $H$-{\it Mumford semi-stable} 
(for $X$ and $H$ as in \S 2) if
$$\mu_H(K) \le \mu_H(E)$$
for all subsheaves $K \subseteq E$ (where $\mu_H$ is the Mumford slope from \S 2). 
A Mumford $H$-semi-stable sheaf $E$ has a {\it Jordan-H\"older filtration}:
$$F_1 \subset F_2 \subset \cdots \subset F_M = E$$
where the $F_{i+1}/F_i$ are Mumford $H$-stable sheaves all of the same slope $\mu_H(E)$.
Although the filtration is not unique, in general, the associated graded coherent sheaf
$\oplus F_i  = \mbox{Ass}_H(E)$ is uniquely determined by the semi-stable sheaf $E$ (and $H$). 

\medskip

The Mumford $H$-slope has the following additional crucial property:

\medskip

\nt {\bf Harder-Narasimhan Filtration:} Every coherent sheaf $E$ on $X$ admits a uniquely determined
(finite) filtration by coherent subsheaves:
$$0 \subset E_0 \subset E_1 \subset E_2 \subset \cdots \subset E_N = E \ \ \ \ \ \mbox{such that}$$

$\bullet$ $E_0$ is the torsion subsheaf of $E$ and

\medskip

$\bullet$ Each $E_{i}/E_{i-1}$ is $H$-semi-stable of slope $\mu_i$ with 
$\mu_1 > \mu_2 > \cdots > \mu_{N}$.

\medskip

Harder-Narasimhan filtrations for a fixed ample divisor class $H$ 
give rise to a family of  ``torsion pairs'' in the category 
of coherent sheaves on $X$:

\medskip

\nt {\bf Definition:} A pair $(\F,\T)$ of full subcategories of a fixed abelian category $\A$ 
is a {\it torsion pair} if:

\medskip

(a) For all objects $T \in \mbox{ob}(\T)$ and $F \in \mbox{ob}(\F)$,  $\mbox{Hom}(T,F) = 0$.

\medskip

(b) Each $A \in \mbox{ob}(\A)$ fits into a (unique) extension
$0 \rightarrow T \rightarrow A \rightarrow F \rightarrow 0$
for some (unique up to isomorphism) objects $T \in \mbox{ob}(\T)$ and $F \in \mbox{ob}(\F)$.

\medskip

\nt {\bf Application 3.1:} For each real number $s$, let $\T_s$ and $\F_s$ be full subcategories of the category 
$\A$ of coherent sheaves on $X$ 
that are closed under extensions and which are generated by, respectively:

\medskip

$\F_s \supset \{ \mbox{torsion-free $H$-stable sheaves of $H$-slope $\mu \le s$}\}$

\medskip

$\T_s \supset \{ \mbox{torsion-free $H$-stable sheaves of $H$-slope $\mu > s$}\} \cup \{\mbox{torsion sheaves}\}$

\medskip

\nt Then $(\F_s,\T_s)$ is a torsion pair of $\A$.

\medskip

{\bf Proof:} Part (a) of the definition follows from the fact that Hom$(T,F) = 0$ if $T, F$ are $H$-stable and $\mu_H(T) > \mu_H(F)$, together with the fact that Hom$(T,F) = 0$ if $T$ is torsion and $F$ is torsion-free. 

\medskip

A coherent sheaf $E$ is either torsion (hence in $\T_s$ for all $s$) or else 
let $E(s) := E_i$ 
be the largest subsheaf in the Harder-Narasimhan of $E$ with the property that $\mu(E_i/E_{i-1}) > s$. Then $0 \rightarrow E(s) \rightarrow E \rightarrow E/E(s) \rightarrow 0$
is the desired short exact sequence for (b) of the definition.

\medskip

\nt {\bf Theorem (Happel-Reiten-Smal\o) \cite{HRS96}:} Given a torsion pair $(\T,\F)$, then there 
is a $t$-structure on the bounded derived category $\D(\A)$ defined by:

\medskip

$\mbox{ob}(\D^{\ge 0}) = 
\{E^\bullet \in \mbox{ob}(\D) \ | 
\mbox{H}^{-1}(E^\bullet) \in \F, \mbox{H}^i(E^\bullet) = 0 \ \mbox{for}\ i < -1 \}$

\medskip

$\mbox{ob}(\D^{\le 0}) = 
\{ E^\bullet \in \mbox{ob}(\D) \ | 
\mbox{H}^0(E^\bullet) \in \T, \mbox{H}^i(E^\bullet) = 0 \ \mbox{for}\ i > 0 \}$

\medskip

\nt 
In particular, the heart of the $t$-structure:  
$$\A_{(\F,\T)} := \{ E^\bullet \ | \ \mbox{H}^{-1}(E^\bullet) \in \F, \mbox{H}^0(E^\bullet) \in \T, 
\mbox{H}^i(E^\bullet) = 0 \ \mbox{otherwise}\}$$
is an  abelian category (referred to as the ``tilt'' of $\A$ with respect to $(\F,\T)$). 

\medskip

\nt {\bf Notation:} We will let $\A_s$ denote the tilt with respect to $(\F_s,\T_s)$ (and $H$).

\medskip

In practical terms, the category $\A_s$ consists of:

\medskip

$\bullet$ Extensions of torsion and $H$-stable sheaves $T$ of slope $> s$

\medskip

$\bullet$ Extensions of shifts $F[1]$ of $H$-stable sheaves $F$ of slope $\le s$ 

\medskip

$\bullet$ Extensions of a sheaf $T$ by a shifted sheaf $F[1]$. 

\medskip

Extensions in $\A_s$ of coherent sheaves $T_1,T_2$ in $\T$ or shifts of coherent sheaves 
$F_1[1],F_2[1]$ in $\F$ are given by extension classes in $\Ext^1_{\A}(T_1,T_2)$ or 
$\Ext^1_{\A}(F_1,F_2)$, which are first extension classes in the category of coherent sheaves.

\medskip

However, an extension of a coherent sheaf $T$ by a shift $F[1]$ in $\A_s$
is quite different. It is given by an element of $\Ext^1_{\A_s}(T,F[1])$ by definition, 
but:
$$\Ext^1_{\A_s}(T,F[1]) = \Ext^2_{\A}(T,F)$$
and this observation will allow us to associate objects of $\A_s$ to certain ``higher'' extension classes
of coherent sheaves in $\A$
just as coherent sheaves are associated to first extension classes of coherent sheaves.

\medskip

First, though, recall that the {\bf rank} is an integer-valued  linear function: 
$$r: K(\D) \rightarrow \bZ$$
on the Grothendieck group of the derived category of coherent sheaves, 
with the property that  $r(E) \ge 0$ for all coherent sheaves $E$.

\medskip

We may define an analogous rank function for each $s \in \bR$ (and $H$):
$$r_s:K(\D) \rightarrow \bR; \ r_s(E) = c_1(E) \cdot H^{n-1} - s \cdot r(E) H^n$$
which has the property that $r_s(E^\bullet) \ge 0$ for all objects 
$E^\bullet$ of $\A_s$ and $r_s(T) > 0$ for all coherent sheaves in $\T_s$ supported in 
codimension $\le 1$.
This rank is 
evidently rational-valued
if $s\in \bQ$.
 
\medskip

Now consider the  objects $\OX[1]$ and $L\otimes \IZ$ of $\A_s$ for $0 \le s < 1$ ($H = c_1(L)$)
where $Z\subset X$ is any closed subscheme supported in codimension $\ge 2$. 

\medskip

\nt {\bf Sub-objects of $\OX[1]$:} 

\medskip

An exact sequence 
$0 \rightarrow K^\bullet \rightarrow \OX[1] \rightarrow Q^\bullet \rightarrow 0$
of objects of $\A_s$ (for any $s \ge 0$)
induces a long exact sequence of cohomology sheaves:
$$0 \rightarrow \mbox{H}^{-1}(K^\bullet) \rightarrow \OX \rightarrow \mbox{H}^{-1}(Q^\bullet) \stackrel \delta \rightarrow 
\mbox{H}^0(K^\bullet) \rightarrow 0$$
Since $\mbox{H}^{-1}(Q^\bullet)$ is torsion-free and $\delta$ is not a (non-zero) isomorphism, then either:

\medskip

($i$) $\mbox{H}^{-1}(K^\bullet) = \OX$ and $Q^\bullet = 0$, $\delta = 0$ and $K^\bullet = \OX[1]$, or else

\medskip

($ii$) $\mbox{H}^{-1}(K^\bullet) = 0$, and $K = K^\bullet$ is a coherent sheaf
with no torsion subsheaf supported in codimension two. The quotient $Q^\bullet = Q[1]$ is 
the shift of a torsion-free sheaf satisfying:
$$0 \le r_s(Q[1]) = -c_1(Q)\cdot H^{n-1} + s\cdot r(Q)H^n < r_s(\OX[1]) = sH^n$$
hence in particular:
$$s\left(1 - \frac 1{r(Q)}\right)  < \mu_H(Q) \le s$$ 

Moreover, if $E$ is a stable coherent sheaf appearing in the associate graded of a semi-stable 
coherent sheaf in the Harder-Narasimhan filtration 
of $Q$, then the same inequality holds for $\mu_H(E)$ (because the $r_s$ rank is additive).

\medskip

\nt {\bf Sub-objects of $L\otimes \IZ$:} 

\medskip

An exact sequence:
$0 \rightarrow {K'} ^\bullet \rightarrow L\otimes \IZ 
\rightarrow {Q'}^\bullet  \rightarrow 0$ in $\A_s$ (for any $s < 1$) induces a long exact sequence
of cohomology sheaves:
$$0 \rightarrow \mbox{H}^{-1}({Q'}^\bullet) \rightarrow \mbox{H}^0({K'}^\bullet) \rightarrow L\otimes \IZ 
\rightarrow \mbox{H}^0({Q'}^\bullet) \rightarrow 0$$
from which it follows that $K' := {K'}^\bullet$ is a torsion-free coherent sheaf, and either:
\medskip

($i'$) $r(K') = 1$, so that $K' = L\otimes \cI_{Z'}$ and ${Q'}^\bullet = L\otimes \left(\cI_Z/\cI_{Z'}\right)$, or else:

\medskip

($ii'$) $r(K') > 1$ and $H^{-1}({Q'}^\bullet) \ne 0$.

\medskip

\nt In either case, we have the inequality:
$$s < \mu_H(K') \le s + \frac{(1-s)}{r(K')}$$
and the same 
inequality when $K'$ is replaced by any $E'$ appearing in the associated graded of a semi-stable 
coherent sheaf in the Harder-Narasimhan filtration of $K'$.

\medskip

\nt {\bf Corollary 3.2:} The alternatives for a non-zero homorphism: 
$$f \in \mbox{Hom}_{\A_s}(L\otimes \IZ,\OX[1]) = \Ext^1_\A(L\otimes \IZ,\OX) \ \ \mbox{for some fixed}\ 0 <  s < 1$$
are as follows:

\medskip

(a) $f$ is injective, with quotient $Q^\bullet = Q[1]$:
$$0 \rightarrow L\otimes \IZ \stackrel f \rightarrow \OX[1] \rightarrow Q[1] \rightarrow  0$$
which in particular implies that $1/2 = \mu_H(Q) \le s$ and, more generally, that each stable $E$
in the Harder-Narasimhan filtration of $Q$ has Mumford-slope $\mu_H(E) \le s$.

\medskip

(b) $f$ is surjective, with kernel $(K')^\bullet = K'$:
$$0 \rightarrow K' \rightarrow L\otimes \IZ \stackrel f\rightarrow \OX[1] \rightarrow 0$$
which in particular implies that $1/2 = \mu_H(K') > s$ and, more generally, that each stable $E'$ in
the Harder-Narasimhan filtration of $K'$ has Mumford-slope $\mu_H(E') > s$.
 
 \medskip
 
 (c) $f$ is neither injective nor surjective, inducing a long exact sequence:
 $$0 \rightarrow L(-D)\otimes \IW \rightarrow L\otimes \IZ \stackrel f\rightarrow \OX[1] \rightarrow (\OX(D)
 \otimes \IV)[1] \rightarrow 0$$
for some effective divisor $D$ satisfying $D\cdot H^{n-1} \le sH^n$ and $D\cdot H^{n-1} < (1-s)H^n$, 
as well as subschemes $V,W \subset X$ supported in codim $\ge 2$. 

\medskip

{\bf Proof:} Immediate from the considerations above.

\medskip

\nt {\it Example:} At $s = 1/2$, we nearly get the same dichotomy as in $\S2$. Here:
$$f   \ \mbox{is injective in} \ \A_{1/2} \Leftrightarrow Q \ \mbox{is $H$-semistable}$$
so the injectivity (or not) of $f$ is equivalent to the semi-stability (or not) of the coherent 
sheaf $E$ expressed as the corresponding (ordinary) extension.

\medskip

Next, recall that the {\bf degree} is an integer-valued linear function:
$$d:K(\D) \rightarrow \bZ; \ \ d(E) = c_1(E) \cdot H^{n-1}$$
(depending upon $H$) with the property that for all coherent sheaves $E$:
$$r(E) = 0 \Rightarrow\big( d(E) \ge 0 \ \mbox{and $d(E) = 0 \Leftrightarrow$ $E$ is supported in codim $\ge 2$}\big)$$

There is an analogous two-parameter family of degree functions ($s\in \bR, t > 0$):
$$d_{(s,t)}: K(\D) \rightarrow \bR; \ \ d_{(s,t)}(E) = \mbox{ch}_2(E)\cdot H^{n-2} - sc_1(E)\cdot H^{n-1} + 
\left(\frac{s^2 - t^2}2\right)r(E)H^n$$
(i.e. a ray of degree functions for each rank $r_s$). Suppose $E^\bullet$ is an object of
$\A_s$ and
$$r_s(E^\bullet) = c_1(E)\cdot H^{n-1} - s \cdot r(E)H^n = 0$$

Then $E^\bullet$ fits into a (unique) exact sequence (in $\A_s$):
$0 \rightarrow F[1] \rightarrow E^\bullet \rightarrow T \rightarrow 0$
where $T$ is a torsion sheaf supported in codimension $\ge 2$, and $F$ is an $H$-semistable coherent
sheaf with $\mu_H(F) = s$. 

\medskip

\nt {\bf Proposition 3.3:} Suppose $r_s(E^\bullet) = 0$ for an object $E^\bullet$  of $\A_s$. Then
for all $t > 0$,
$$d_{(s,t)}(E^\bullet) \ge 0 \ \mbox{and}\ d_{(s,t)}(E^\bullet) = 0 \Leftrightarrow E^\bullet 
\ \mbox{is a sheaf, supported in codim $\ge 3$}$$

\medskip

{\bf Proof:} Because $d_{(s,t)}$ is linear, it suffices to prove the Proposition for torsion sheaves $T$ supported in 
codimension $\ge 2$ and for shifts $F[1]$ of {\it $H$-stable} torsion-free sheaves  of slope $s$. In the 
former case:
$$d_{(s,t)}(T) = \mbox{ch}_2(T)\cdot  H^{n-2} \ge 0 \ \mbox{with equality $\Leftrightarrow T$ is supported in codim
$\ge 3$}$$

In the latter case:
$$d_{(s,t)}(F[1]) = -\mbox{ch}_2(F)\cdot H^{n-2} + sc_1(F)\cdot H^{n-1} - \left(\frac{s^2-t^2}2\right)r(F)H^n $$
and $\mu_H(F) = s \Rightarrow (c_1(F) - sr(F)H)\cdot H^{n-1} = 0 \Rightarrow 
(c_1(F) - sr(F)H)^2H^{n-2} \le 0$
by the Hodge Index Theorem. 
It follows from the Bogomolov inequality  that:
$$d_{(s,t)}(F[1]) \ge - \left(\frac{c_1^2(F)}{2r(F)}\right)\cdot H^{n-2} + sc_1(F)\cdot H^{n-1} - 
\left(\frac{s^2}2\right)r(F)H^n + \left(\frac{t^2}2\right)r(F)H^n$$
$$\ \ \ \ \ \ \ \ = - \left(\frac 1{2r(F)}\right)(c_1(F) -  sr(F)H)^2\cdot H^{n-2} + \left(\frac{t^2}{2}\right)r(F)H^n > 0 \qed$$

\medskip

\nt {\bf Corollary 3.4:} If $X = S$ is a surface, then the complex linear function:
$$Z_{s+it} := (-d_{(s,t)} + it r_s): K(\D) \rightarrow \bC; \ s\in \bR, t  > 0, i^2 = -1$$
has the property that $Z_{s+it}(E^\bullet) \ne 0$ for all objects $E^\bullet \ne 0$  of $\mbox{ob}(\A_s)$, and:
$$0 < \arg(Z_{s+it}(E^\bullet)) \le 1 \ \mbox{(where $\arg(re^{i\pi \rho}) = \rho$)}$$
i.e. $Z_{s+it}$ takes values in the (extended) upper half plane.

\medskip

In higher dimensions, the Corollary holds modulo coherent sheaves supported in codimension $\ge 3$,
 just as the ordinary $H$-degree and rank lead to the same conclusion modulo torsion sheaves supported in codimension $\ge 2$.

\medskip

\nt {\it Remark:} The ``central charge'' $Z_{s+it}$ has the form:
$$Z_{s+it}(E) = -d_{(s,t)}(E) + itr_s(E) = -\int_S e^{-(s+it)H} \mbox{ch}(E)H^{n-2}$$
which is a much more compact (and important) formulation.

\medskip

\nt {\bf Corollary 3.5:} Each ``slope'' function:
$$\mu := \mu_{s+it}  =  \frac{d_{(s,t)}}{tr_s} = -\frac{\mbox{Re}(Z_{s+it})}{\mbox{Im}(Z_{s+it})}$$ 
has the usual properties 
of a slope function on the objects of  $\A_s$. That is, given an exact sequence of objects
of $\A_s$:
$$0 \rightarrow K^\bullet \rightarrow E^\bullet \rightarrow Q^\bullet \rightarrow 0$$
then
$\mu(K^\bullet) < \mu(E^\bullet) \Leftrightarrow \mu(E^\bullet) < \mu(Q^\bullet)$ and 
$\mu(K^\bullet) = \mu(E^\bullet) \Leftrightarrow \mu(E^\bullet) = \mu(Q^\bullet)$

\medskip

Also, when we make the usual:

\medskip

\nt {\bf Definition:} $E^\bullet \ \mbox{is $\mu$-stable if}\ \mu(K^\bullet) < \mu(E^\bullet) \ \mbox{whenever}
\ K^\bullet \subset E^\bullet$ and the quotient has nonzero central charge
(i.e. is not a torsion sheaf supported in codim $\ge 3$).

\medskip

Then Hom$(E^\bullet,F^\bullet) = 0$ whenever $E^\bullet,F^\bullet$ are $\mu$-stable and 
$\mu(E^\bullet) > \mu(F^\bullet).$

\medskip

{\bf Proof (of the Corollary):} Simple arithmetic. 

\medskip

\nt {\it Example:} In dimension $n \ge 2$:
$$\mu_{s+it}(\OX[1]) = \frac{t^2 - s^2}{2st} \ \ \mbox{and}\ \ \mu_{s+it}(L\otimes \IZ) = \frac{(1-s)^2 - t^2 - 
\frac{2d}{H^n}}{2t(1-s)}$$
where $d = [Z] \cap H^{n-2}$ is the (codimension two) degree of the subscheme $Z \subset X$.
Thus $\mu_{s+it}(L\otimes \IZ) > \mu_{s+it}(\OX[1])$ if and only if:
$$t^2 + \left(s-\left(\frac 12 - \frac d{H^n}\right)\right)^2 < \left(\frac 12 - \frac d{H^n}\right)^2 \ \mbox{and} \ t > 0.$$
This describes a nonempty subset (interior of a semicircle) of $\bR^2$  if $H^n > 2d$. 

\medskip

\nt {\bf Proposition 3.6:} For all smooth projective varieties $X$ of dimension $\ge 2$ (and $L$)

\medskip

(a) $\OX[1]$ is a $\mu_{s+it}$-stable object of $\A_s$ for all $s \ge 0 $ and $t>0$.

\medskip

(b) $L$ is a $\mu_{s+it}$-stable object of $\A_s$ for all $s < 1$ and $t > 0$.

\medskip

{\bf Proof:} (a) Suppose $0 \ne K^\bullet \subset \OX[1]$, 
and let $E$ be an $H$-stable torsion-free sheaf in the associated graded of $Q$, where $Q[1]$ is the 
quotient object. Recall that
$0 < \mu_H(E) \le s$.
The Proposition follows once we show $\mu_{s+it}(\OX[1]) < \mu_{s+it}(E[1])$ for all 
$E$ with these properties.
We compute:
$$\mu_{s+it}(E[1]) = \frac{-2\mbox{ch}_2(E)H^{n-2} + 2sc_1(E)H^{n-1} - \left({s^2-t^2}\right)r(E)H^n}
{2t(-c_1(E)H^{n-1} + sr(E)H^n)}$$
and we conclude (using the computation of $\mu_{s+it}(\OX[1])$ above) that:
$$\mu_{s+it}(\OX[1]) > \mu_{s+it}(E[1]) \Leftrightarrow 
(s^2+t^2)c_1(E)H^{n-1} > (2s)\mbox{ch}_2(E)H^{n-2}$$
But by the Bogomolov Inequality:
$$(2s)\mbox{ch}_2(E)H^{n-2} \le s (c_1^2(E)H^{n-2})/r(E)$$
 and by the Hodge Index Theorem and the inequality $c_1(E)\cdot H^{n-1} \le sr(E)H^n$:
$$sc_1^2(E)H^{n-2}/r(E) \le s^2c_1(E)H^{n-1}$$
The desired inequality follows from the fact that $t > 0$ and $c_1(E)\cdot H^{n-1} > 0$. 

\medskip

The proof of (b) proceeds similarly. Suppose $0 \ne (K')^\bullet \subset L$ in $\A_s$, and let $E'$ be an
$H$-stable coherent sheaf in the Harder-Narasimhan filtration of $K'$. 

Then $s < \mu_H(E') < 1$, and
we need to prove that $\mu_{s+it}(E') < \mu_{s+it}(L)$. 
This follows as in (a) from the Bogomolov Inequality and Hodge Index Theorem. $\qed$.

\medskip

\nt {\bf Corollary 3.7:} (Special case of Kodaira vanishing):
$$\mbox{H}^{n-1}(X,K_X+L) = 0 \ \mbox{for all $n > 1$} $$

{\bf Proof:} Within the semicircle
$\{ (s,t) \ | \ t^2 + \left(s-\frac 12\right)^2 < \frac 14 \ \mbox{and}\ t > 0\}$
the inequality $\mu_{s+it}(L) > \mu_{s+it}(\OX[1])$ holds. But $L$ and $\OX[1]$ are always $\mu_{s+it}$-stable, hence: 
$$0 = \mbox{Hom}_{\A_s}(L,\OX[1]) \cong \Ext^1_{\OX}(L,\OX) \cong \mbox{H}^{n-1}(X,K_X+L)^\vee$$
\qed

\medskip

\nt {\it Remark:} The Bogomolov Inequality and Hodge Index Theorem are trivially true in dimension one.
However, the computation of $\mu_{s+it}(L)$ is different in dimension one, and indeed in that case the inequality $\mu_{s+it}(L) > \mu_{s+it}(\OX[1])$ never holds (which is good, since the corollary is false in dimension one)!

\medskip

Restrict attention to $X = S$ a surface for the rest of this section, and consider:
$$\IW^\vee[1]$$
the shifted derived dual of the ideal sheaf of a subscheme $W \subset S$ of length $d$. Since:
$$\mbox{H}^{-1}(\IW^\vee[1]) = \OS \ \ \mbox{and}\ \ \mbox{H}^0(\IW^\vee[1]) \ \mbox{is a torsion sheaf, supported on $W$}$$
it follows that $\IW^\vee[1]$ is in $\A_s$ for all $s \ge 0$. 

\medskip

Every quotient object $\IW^\vee[1] \rightarrow Q^\bullet$ satisfies:

\medskip

$\bullet$ $H^{0}(Q^\bullet)$ is supported in codimension two (on the scheme $W$, in fact).

\medskip

$\bullet$ Let $Q = H^{-1}(Q^\bullet)$ (a torsion-free sheaf). Then every $H$-stable term $E$ in the Harder-Narasimhan filtration of $Q$ satisfies:
$$0 \le r_s(E[1]) =  -c_1H + rsH^2 < r_s(\IW^\vee[1]) = sH^2 \Leftrightarrow (r-1)sH^2 < c_1H \le rsH^2$$
where $r = r(E)$ and  $c_1 = c_1(E)$ (because the $r_s$ rank of the kernel object is positive).

\medskip

\nt {\bf Proposition 3.8:} For subschemes $Z,W \subset S$ of the same length $d$ (and $H = c_1(L)$):

\medskip

(a) If $H^2 > 8d$, then $\mu_{s+it}(L\otimes\IZ) > \mu_{s+it}(\IW^\vee[1])$ for all $(s,t)$ in the semicircle:
$$C(d,H^2) := \left\{(s,t) \ | \ t^2 + \left(s-\frac 12\right)^2 < \frac 14 - \frac {2d}{H^2} \ \mbox{and}\ t > 0\right\}$$ 
centered at the point $(1/2,0)$ (and the semicircle is nonempty!).

\medskip

(b)  If  $H^2 > 8d$ and $\IW^\vee[1]$ or $L\otimes \IZ$ is not stable at 
$(s,t) = (\frac 12, \sqrt{\frac 14 - \frac {2d}{H^2}})$, then there is a divisor $D$ 
on $S$ and an integer $r > 0$ such that:
$$\frac{r-1}2H^2 < D\cdot H \le \frac r2H^2, \ \ \mbox{and}\ \ \frac Dr \cdot H < \frac{D^2}{r^2} + 2d$$

{\bf Proof:} Part (a) is immediate from:
$$\mu_{s+it}(L\otimes \IZ) = \frac{(1-s)^2 - t^2 - \frac{2d}{H^2}}{2(1-s)t} \ \mbox{and}\ \mu_{s+it}(\IW^\vee[1]) = \frac{t^2 - s^2 + \frac{2d}{H^2} }{2st}$$

We prove part (b) for $\IW^\vee[1]$
(the proof for $L\otimes \IZ$ is analogous). 

\medskip

Let $\IW^\vee[1] \rightarrow Q^\bullet$ be a surjective map in the category $\A_s$ and let $Q = H^{-1}(Q^\bullet)$. Since 
$H^0(Q^\bullet)$ is torsion, supported on $W$, it follows that
$\mu_{s+it}(Q[1]) \le \mu_{s+it}(Q^\bullet)$ with equality if and only if $H^0(Q^\bullet) = 0$. 

\medskip

Thus  
if $\IW^\vee[1]$ is not $\mu_{s+it}$-stable, then $\mu_{s+it}(\IW^\vee[1]) \ge \mu_{s+it}(Q[1])$ for some 
torsion-free sheaf $Q$ satisfying $(r-1)sH^2 < c_1(Q)\cdot H \le rsH^2$, and moreover, the same set of 
inequalities hold for (at least) one of the stable torsion-free sheaves $E$ appearing in the Harder-Narasimhan filtration of $Q$.
We let $D = c_1(E)$ and $r = \mbox{rk}(E)$. Then
$\mu_{s+it}(\IW^\vee[1]) \ge \mu_{s+it}(E[1])$ if and only if:
$$(t^2 + s^2)(D\cdot H) \le (2s)\mbox{ch}_2(E) + \frac{2d}{H^2} (rsH^2 - D\cdot H)$$
and by the Bogomolov inequality, $(2s)\mbox{ch}_2(E) \le s\frac{D^2}r$. Setting $(s,t) = (\frac 12, \sqrt{\frac 14 - \frac {2d}{H^2}})$, we obtain the desired inequalities. \qed

\medskip
 
\nt {\bf Corollary 3.9:} (a) If $L = \OS(H)$ is ample on $S$ and satisfies $H^2 > (2d+1)^2$ and:
$$\mbox{H}^1(S,K_S\otimes L \otimes I_W\otimes I_Z) \ne 0$$
for a pair $Z,W \subset S$ of length $d$ subschemes then there 
is a divisor $D$ on $S$ satisfying $D^2 \le 0$ and $0 < D\cdot H \le D^2 + 2d$.

\medskip

(b) (Fujita-type result) If $L$ is an arbitrary ample line bundle on $S$, then
$$\mbox{H}^1(S,K_S\otimes L^{\otimes (2d+2)} \otimes I_W\otimes I_Z) = 0$$
for all subschemes $Z,W \subset S$ of length $d$ (or less).

\medskip

{\bf Proof:} Part (b) immediately follows from (a). Since: 
$$\mbox{H}^1(S,K_S\otimes L \otimes I_W\otimes I_Z) \cong 
\mbox{Hom}_{\A_{\frac 12}}(L\otimes \IZ, \IW^\vee[1])^\vee$$
and $(2d+1)^2 \ge 8d + 1$ for all $d \ge 1$, the non-vanishing of $H^1$ implies 
that either $L\otimes \IZ$ or $\IW^\vee[1]$ must not be 
stable at $(\frac 12, \sqrt{\frac 14 - \frac {2d}{H^2}})$, and from Proposition 3.8 (b) there is a divisor
$D$ and integer $r \ge 1$ such that the $\bQ$-divisor $C = D/r$ satisfies:
$$(1 - \frac 1r)\frac{H^2}2 < C\cdot H \le \frac{H^2}2 \ \ \mbox{and}\ \ C\cdot H \le C^2 + 2d$$
(similar to Application 2.2). The result now follows as in Application 2.2 once we prove that
$C^2 \ge 1$ whenever $r > 1.\footnote{The authors thank Valery Alexeev for pointing out the embarassing 
 omission of this step in the original version of the paper.}$ To this end, note:
 
 \medskip
 
(i) $r \ge 3 \Rightarrow C^2 + 2d \ge C\cdot H > \frac{H^2}{3} > \frac{8d + 1}3 \Rightarrow C^2 >
\frac{2d + 1}3  \ge 1$.

\medskip

(ii) $r = 2 \Rightarrow C^2 + 2d \ge C \cdot H > \frac{H^2}4 > 2d + \frac 14 \Rightarrow \frac{H^2}4 \ge 2d + \frac 12, C\cdot H \ge 2d + 1$, and $C^2 \ge 1$, since $C$ is of the form $D/2$ for an ``honest''
divisor $D$.

\medskip

Thus either $C^2 \le 0$, in which case $r = 1$ and $C = D$ is an ``honest'' divisor, or else $C^2 \ge 1$. Furthermore, by the Hodge index theorem:
$$C^2H^2 \le (C\cdot H)^2  \le \frac{H^2}2\left(C^2 + 2d\right) \Rightarrow C^2 \le 2d$$
and if $C^2 = \kappa$ for $1 \le \kappa \le 2d$, then
$\kappa^2H^2 \le (C\cdot H)^2 \le (\kappa + 2d)^2 \Rightarrow H^2 \le \left(1 + \frac{2d}{\kappa}\right)^2$.
This is a decreasing function, giving us $H^2 \le (2d+1)^2$, contradicting
$H^2 > (2d+1)^2$. 

\medskip

\nt {\it Remark:} This variation resembles other variations of Reider's theorem, e.g.  \cite{Lan99}, though 
the authors do not see how to directly obtain this result from the others.

\medskip

In a special case, Proposition 3.8 can be made even stronger, as noted in \cite{AB}.

\medskip

\nt {\bf Proposition 3.10:} If Pic$(S) = \bZ$, generated by $c_1(L) = H$, then the two objects
$L\otimes \IZ$ and $\IW^\vee[1]$ are 
$\mu_{(\frac 12,t)}$-stable for all $t > 0$ and any degree of $Z$ (and $W$).

\medskip

{\bf Proof:} Again we do this for $\IW^\vee[1]$, the proof for $L\otimes \IZ$ being analogous. Consider again 
the condition on every subbundle $E \subset Q$, where $Q = H^{-1}(Q^\bullet)$, 
and $Q^\bullet$ is a quotient object of $\IW^\vee[1]$:
$$(r(E)-1)(\frac 12)H^2 < c_1(E)\cdot H \le r(E)(\frac 12)H^2$$
Since $c_1(E) = kH$ is an integer multiple of $H$, by assumption, it follows immediately that $Q$ is itself 
of even rank and $H$-stable, satisfying $c_1(Q) = \left(r(Q)/2\right)H$. But in that case,
$Q[1]$ has ``Bridgeland rank'' $r_{\frac 12}(Q[1]) = 0$, hence has maximal phase (infinite slope), and thus cannot destabilize $\IW^\vee[1]$. \qed

\medskip

\nt {\it Remark:} This argument is highly sensitive to setting $s = \frac 12$, and indeed the conclusion is
not true when $s \ne \frac 12$.

\medskip

\nt {\bf Corollary 3.11:} If Pic$(S) = \bZ H$ and $H^2 > 8d$, then $H^1(S,K_S\otimes L \otimes \IW\otimes \IZ) = 0$
for pairs of subschemes $Z,W \subset S$ of length $d$.

\medskip

\nt {\bf 4. Thaddeus Pairs Revisited.} Let $S$ be a surface with ample line bundle $L$ and Pic$(S) = \bZ H$ 
with $H = c_1(L)$.
Consider the objects of $\A_s$ ($0 < s < 1$) appearing as extensions:
$$\epsilon: 0 \rightarrow \cO_S[1] \rightarrow E^\bullet_\epsilon \rightarrow L \rightarrow 0$$
parametrized by:
$$\epsilon \in \mbox{Ext}^1_{\A_s}(L,\cO_S[1]) = 
\mbox{Ext}^2_{\cO_S}(L,\cO_S) \cong \mbox{H}^0(S,K_S\otimes L)^\vee$$

As we saw in Proposition 3.6 and the preceding calculation, $\cO_S[1]$ and $L$ are both $\mu_{s+it}$
stable for all $(s,t)$. Moreover, 
$\mu_{s+it}(\cO_S[1]) < \mu_{s+it}(L)$ inside the semicircle:
$$C := \left\{ (s,t) \ | \ t^2 + \left(s-\frac 12\right)^2 < \frac 14 \ \mbox{and} \ t > 0\right\}$$

\nt {\it Remark:} Here and earlier, we are using the notion of stability a little bit loosely. The correct definition, given by Bridgeland \cite{Bri08} requires the existence of finite-length Harder-Narasimhan filtrations for all objects of 
$\A_s$. This is straightforward to prove when $(s,t)$ are both rational numbers (following Bridgeland), but much more subtle in the irrational case. For the purposes of this paper, the rational values will suffice.

\medskip

We investigate the dependence of the $\mu_{\frac 12+it}$-stability of $E^\bullet_\epsilon$ upon the extension class $\epsilon$ for $\frac 12 + it$ inside the semicircle $S$. If $E^\bullet_\epsilon$ is $\mu_{\frac 12 + it}$-unstable, destabilized by 
$$K^\bullet \subset E^\bullet_\epsilon, \ \mbox{then:}$$

\medskip

(i) $K^\bullet = \mbox{H}^0(K^\bullet) =: K$ is a coherent sheaf with $\mu_{\frac 12+it}(K) > 0$.

\medskip

(ii) $K$ is $H$-stable of odd rank $r$ and $c_1(K) = ((r+1)/2)H$.

\medskip

(iii) The induced map $K \rightarrow L$ is injective (in the category $\A_{\frac 12}$).

\medskip

Thus as in the curve case, $E^\bullet_\epsilon$ can only be destabilized by lifting subobjects $K \subset L$ 
(in the category $\A_{\frac 12}$) of 
positive $\mu_{\frac 12 + it}$-slope to subobjects of $E^\bullet_\epsilon$:
$$\begin{array}{ccccccccc}
&&&&&& K\\
&&&&&\swarrow & \downarrow\\
(\dag) \ \ \  0 & \rightarrow & \OS[1] & \rightarrow & E^\bullet_\epsilon & \rightarrow & L & \rightarrow & 0
\end{array}$$
That is, the unstable objects $E^\bullet_\epsilon$ correspond to extensions in the kernel of the map:
$$\mbox{Ext}^2(L,\OS) \rightarrow \mbox{Ext}^2(K,\OS)$$
for some mapping of coherent sheaves $K \rightarrow L$ with $K$ satisfying (i) and (ii).

\medskip

{\bf Proof} (of (i)-(iii)): The $d_{(\frac 12,t)}$-degree of $E^\bullet_\epsilon$ is:
$$\mbox{ch}_2(E^\bullet_\epsilon) - \frac {c_1\cdot H}2 + \frac {(\frac 14 - t^2)}2 rH^2 = 0$$
since $\mbox{ch}_2(E^\bullet_\epsilon) = H^2/2, c_1 = H$ and $r=0$. Thus the slope (equivalently, the degree) of any destabilizing $K^\bullet \subset E^\bullet_\epsilon$ is positive, by definition. Moreover the ``ranks''
$$r_{\frac 12}(\OS[1]) = r_{\frac 12}(L) = \frac {H^2}{2}$$
are the minimal possible (as in the curve case) without being zero, hence as in the curve case, 
$K^\bullet \subset E^\bullet_\epsilon$ must also have minimal rank $\frac {H^2}2$ (if it had the next smallest rank 
$H^2 = r_{\frac 12}(E^\bullet_\epsilon)$, it would fail to destabilize). The presentation of $E^\bullet_\epsilon$ gives 
$H^{-1}(E^\bullet_\epsilon) = \OS$ and $H^0(E^\bullet_\epsilon) = L$, hence if we let $Q^\bullet = E^\bullet_\epsilon/K^\bullet$, then:
$$0 \rightarrow H^{-1}(K^\bullet) \rightarrow \cO_S \rightarrow H^{-1}(Q^\bullet) \rightarrow
K \rightarrow L \rightarrow H^0(Q^\bullet) \rightarrow 0$$
and, as usual, either $H^{-1}(K^\bullet) = 0$ or $H^{-1}(K^\bullet) = \OS$. The latter is impossible, since
in that case, the rank consideration would give $K^\bullet = \OS[1]$, which doesn't destabilize for $(\frac 12,t) \in C$. Thus 
$K^\bullet = K$ is a coherent sheaf. This gives (i).

\medskip

Next, the condition that $r_{\frac 12}(K)$ be minimal implies that there can 
only be one term in the Harder-Narasimhan filtration of $K$ (i.e. $K$ is $H$-stable), and 
that:
$$r_{\frac 12}(K) = c_1(K)H - \frac  {r(K)H^2}2 = \frac {H^2}2.$$
Since $c_1(K) = kH$ for some $k$, this gives (ii).
 
 \medskip
 
Finally, (iii) follows again from the minimal rank condition since any kernel of the induced map to $L$ would be a torsion-free sheaf,
of positive $r_{\frac 12}$-rank. \qed

\medskip

Suppose now that $K$ satisfies (i) and (ii). By the Bogomolov inequality:
$$d_{(\frac 12,t)}(K)  \le \frac 1{2r}\left(c_1(K) - \frac r2H\right)^2 - \frac {rt^2H^2}2 = \frac {H^2}{2r}\left(\frac 14 - r^2t^2\right)$$
so in particular, $t \le \frac 1{2r}$, or in other words, we have shown:

\medskip

\nt {\bf Proposition 4.1:} If $t > \frac 1{2r}$ and $\mu_{\frac 12 + it}(K) < 0$ for all $K \subset L$ (in $\A_{\frac 12}$)
of odd ordinary rank $ \le r$, then $E^\bullet_\epsilon$ is $\mu_{\frac 12 + it}$-stable.

\medskip

\nt {\bf Special Case:} Suppose $t > \frac 16$. Because $H = c_1(L)$ generates Pic$(S)$ it follows that the 
only rank one subobjects $K \subset L$ in $\A_{\frac 12}$ are the subsheaves $L\otimes \IZ$ for $Z \subset S$
of finite length. Thus
$E^\bullet_\epsilon$ only fails to be $\mu_{\frac 12 + it}$-stable if:
$$d_{(\frac 12,t)}(L\otimes \IZ) = \frac 12\left(\frac 14 - t^2\right)H^2 - d \ge 0 \Leftrightarrow t^2 \le \frac 14 - \frac{2d}{H^2}$$
and $\epsilon \in \mbox{ker}(\mbox{Ext}^2(L,\OS) \rightarrow \mbox{Ext}^2(L\otimes \IZ,\OS))$, so that 
$L\otimes \IZ \subset L$ lifts to a subobject of $E^\bullet_\epsilon$. As in the curve case,
Serre duality implies that the image of such a (non-zero) extension in the projective space:
$$\bP(H^0(S,K\otimes L)^\vee)$$
is a point of the {\it secant $d-1$-plane} spanned by $Z \subset S$ under the linear series map:
$$\phi_{K+L}: S --> \bP(H^0(S,K\otimes L)^\vee).$$

By Corollary 3.11, this inequality on $t$ guarantees 
$\mbox{H}^1(S,K_S\otimes L \otimes \IW\otimes \IZ)^\vee = 0$
for all subschemes $Z,W \subset S$ of length $d$, hence in particular, the $d-1$-secant planes spanned by $Z \subset S$ are well-defined.

\medskip

Thus there are ``critical points''  or ``walls'' at
$t = \sqrt{\frac 14 - \frac{2d}{H^2}} > \frac 16, \ \mbox{i.e.}\ d < \frac {2H^2}9$
on the line $s = \frac 12$ where the objects $E^\bullet_\epsilon$ corresponding to points of the secant variety:
$$\left(\mbox{Sec}^{d-1}(S) - \mbox{Sec}^{d-2}(S)\right) \subset \bP(H^0(S,K\otimes L)^\vee)$$
change from $\mu$-stable to $\mu$-unstable as $t$ crosses the wall.

\medskip

\nt {\bf Moduli.} The Chern class invariants of each $E^\bullet_\epsilon$ are:
$$\mbox{ch}_2 = \frac {H^2}2, \ \ \ c_1 = H, \ \ r = 0$$
Thus it is natural to ask for the set of all $\mu_{\frac 12 + it}$-stable objects with these invariants, and 
further to ask whether they have (projective) moduli that are closely related (by flips or flops) as 
$t$ crosses over a critical point. In one case, this is clear:

\medskip

\nt {\bf Proposition 4.2:} For $t > \frac 12$, the $\mu_{\frac 12 +it}$-stable objects with Chern class invariants 
above are precisely the (Simpson)-stable coherent sheaves with these invariants, i.e. sheaves of 
pure dimension one and rank one on curves in the linear series $|H|$.

\medskip

{\bf Proof:} Suppose $E^\bullet$ has the given invariants and is not a coherent sheaf. Then the sequence:
$$0 \rightarrow \mbox{H}^{-1}(E^\bullet)[1] \rightarrow E^\bullet \rightarrow \mbox{H}^0(E^\bullet) \rightarrow 0$$
destabilizes $E^\bullet$ for $t > \frac 12$ for the following reason. Let $E = H^{-1}(E^\bullet)$. If $c_1(E) = kH$, then $k \le \frac r2$ is required in order that $E[1] \in \A_{\frac 12}$. Moreover, since $r_{\frac 12}(E^\bullet) = H^2$
and $H^0(E^\bullet)$ has positive (ordinary) rank, hence also positive $r_{\frac 12}$-rank, it follows that 
$r_{\frac 12}(E[1]) =  0 \ \mbox{or} \ \frac{H^2}2$. But $r_{\frac 12}(E[1]) = 0$ implies $E[1]$ has maximal
(infinite) slope, and then $E^\bullet$ is unstable (for all $t$). It follows similarly that if $r_{\frac 12}(E[1]) = \frac{H^2}2$, then $E^\bullet$ is unstable for all $t$ unless $E$ is $H$-stable, of rank $r = 2k + 1$. 
In that case, by the Bogomolov inequality:
$$d_{(\frac 12,t)}(E[1]) = -\mbox{ch}_2(E) + \frac {c_1(E)H}2 - \frac{(\frac 14 - t^2)rH^2}2
\ge \frac {H^2}2 \left(t^2r - \frac 1{4r}\right)$$
and this is positive if $t > \frac 12$. \qed

\medskip

In fact, at $t = \frac 12$ only $H^{-1}(E^\bullet) = \OS$ (the rank one case, moreover matching the Bogomolov bound) would fail to destabilize a non-sheaf $E^\bullet$, and 
conversely, among the coherent sheaves $T$ with these invariants, only those fitting into an exact sequence
(of objects of $\A_{\frac 12})$:
$$0 \rightarrow L \rightarrow T \rightarrow \OS[1] \rightarrow 0$$
become unstable as $t$ crosses $\frac 12$, and they are replaced by the ``Thaddeus'' extensions:
$$0 \rightarrow \OS[1] \rightarrow E^\bullet \rightarrow L\rightarrow 0$$

In other words, the moduli of Simpson-stable coherent sheaves
$$M_S\left(0,H,\frac{H^2}2\right)$$
is known to be projective by a geometric invariant theory construction \cite{Sim94}. 
It is the moduli of $\mu_{\frac 12 + it}$-stable objects of $\A_{\frac 12}$
for $t > \frac 12$. The wall crossing at $t = \frac 12$ removes:
$$\bP(\mbox{Ext}^1(\OS[1],L)) = \bP(H^0(S,L)) \subset M_S\left(0,H,\frac{H^2}2\right)$$
and replaces it with $\bP(H^0(S,K_S+L)^\vee)$ in another birational model. In the case $K_S = 0$,
the Simpson moduli spaces are holomorphic symplectic varieties, this new birational model is a Mukai flop
of the moduli of stable sheaves, and the further wall crossings (up to $t = \frac 16$, when rank three bundles appear) all replace extensions of the form:
$$0 \rightarrow L\otimes \IZ \rightarrow (T \ \mbox{or $E^\bullet$}) \rightarrow \IW^\vee[1] \rightarrow 0$$
with
$$0 \rightarrow  \IW^\vee[1] \rightarrow E^\bullet \rightarrow L\otimes \IZ \rightarrow 0.$$

This is achieved globally by Mukai flops, replacing projective 
bundles over the product $\mbox{Hilb}^d(S) \times \mbox{Hilb}^d(S)$ of Hilbert schemes with their dual bundles:
$$\bP(H^0(S,L\otimes \IW\otimes \IZ)) \ \leftrightarrow \bP(H^0(S,L\otimes \IW\otimes \IZ)^\vee)$$
This was constructed in detail in \cite{AB}.

\medskip

General questions regarding moduli of Bridgeland-stable objects remain fairly wide open, however.
Toda \cite{Tod08} has shown that when $S$ is a $K3$ surface, then the Bridgeland semistable objects 
of fixed numerical class are represented by an Artin stack of finite type. One expects the isomorphism classes of Bridgeland-stable objects, at least in special cases as above, to be represented by a proper scheme when 
$(s,t)$ is not on a ``wall.'' However:

\medskip

\nt {\bf Question 1:} When are the isomorphism classes of Bridgeland-stable objects of fixed numerical type represented by a (quasi)-projective scheme of finite type?

\medskip

\nt {\bf Question 2:} Conversely, is there an example where the isomorphism classes are represented by 
a proper algebraic space which is not a projective scheme? (The examples produced in \cite{AB} are proper algebraic spaces. It is unknown whether they are projective.)

\medskip

For each $t < \frac 12$, we make the following provisional:

\newpage

\nt {\bf Definition:} The space of $t$-stable Thaddeus pairs 
(given $S$ and ample $L$) is the proper transform of 
the projective space of extensions $\bP(\mbox{Ext}^1(L,\OS[1]))$ under the natural rational embedding in the moduli space of (isomorphism classes of) $\mu_{\frac 12 +it}$-stable 
objects with invariants $(0,H,H^2/2)$.

\medskip

\nt {\it Remark:} Note that for $t < \frac 16$, this will contain objects that have no analogue in the curve 
case, corresponding to destabilizing Mumford-stable torsion-free sheaves $K$ of higher odd rank $r$ and
first chern class $c_1(K) = \frac{r+1}2H$. 

\medskip

\nt {\bf Question 3.} Can stable Thaddeus pairs, as a function of $t$ (inside the moduli of $\mu_{\frac 12 + it}$-stable objects of the same numerical class) be defined as a moduli problem? If so, what are its properties? 
Is it projective? Smooth? What happens as $t \downarrow 0$?

\bigskip

\nt {\bf 5. Reider in Dimension Three?} Let $X$ be a smooth projective three-fold, with Pic$(X) = \bZ \cdot H$
(for simplicity) and $H = c_1(L)$ for an ample line bundle $L$. Consider:
$$\epsilon: 0 \rightarrow \OX[1] \rightarrow E^\bullet_\epsilon \rightarrow L\otimes \IZ \rightarrow 0$$
for subschemes $Z \subset X$ of finite length $d$, taken within the tilted category $\A_{\frac 12}$.

\medskip

\nt {\bf Question 4:} Are there bounds $d_0$ and $t_0$ such that {\it all} objects 
$E^\bullet_\epsilon$ formed in this way are $\mu_{\frac 12 + it}$-unstable when $d > d_0$
and $t < t_0$? If so, does this follow from a more general Bogomolov-type codimension three
inequality for the numerical invariants of
$\mu$-stable objects?

\medskip

As we have already discussed in the surface case, a destabilizing subobject of an  {\it unstable} such $E^\bullet_\epsilon$ would be exhibited by lifting:
$$\begin{array}{ccccccccc}
&&&&&& K\\
&&&&&\swarrow & \downarrow\\
(\dag) \ \ \  0 & \rightarrow & \OS[1] & \rightarrow & E^\bullet_\epsilon & \rightarrow & L\otimes \IZ & \rightarrow & 0
\end{array}$$
where $K \subset L\otimes \IZ$ is a subobject in $\A_{\frac 12}$. By requiring Pic$(X) = \bZ$, we assure that 
$K$ does not factor through $L(-D)$ for any effective divisor $D$. The interesting cases are therefore:

\medskip

$\bullet$ $K = L\otimes \IC$ where $C \subset X$ is a curve, and a new-looking condition:

\medskip

$\bullet$ $K$ is an $H$-stable torsion-free sheaf of odd rank $r$ and $c_1(K) = \frac{r+1}2H$.

\medskip

\nt {\bf Question 5:} Assuming Question 4, are there examples of threefolds where the bounds of 
Question 4 
are satisfied (hence all $E^\bullet$ are $\mu$-unstable), but the ``interesting cases'' above allow for non-zero extensions?

\medskip

And the last question is whether all the interesting cases can be numerically eliminated (by some form of the Hodge Index Theorem) leading to a proof of Fujita's conjecture for threefolds.

\end{document}